\documentclass[12pt,a4paper]{article}
 \usepackage{emlines2,amstext,amssymb,amscd}

\begin{document}

\title{A simple proof of associativity and commutativity of LR-coefficients (or the hive ring)}
\author{    V.Danilov\thanks{The support of the grant NSh-1939.2003.6, School Support, is grateful
acknowledged.} \ and G.Koshevoy\thanks{ The support of LIFR MIIP
is grateful
acknowledged.}\\ Central Institute of Economics and Mathematics\\
Nahimovskii pr. 47, 117418 Moscow, Russia\\ e-mail:
vdanilov43@mail.ru, koshevoy@cemi.rssi.ru }

\maketitle

In this paper we propose a simple bijective proof of associativity
and commutativity of Littlewood-Richardson coefficient or the hive
ring (\cite{KTW2}).

\section{Introduction}

The ring of symmetric functions $\mathbb Z [x_1,\ldots, x_n]$ has
a distinguished basis constituted of  the Schur functions
$s_{\lambda }$, where  $\lambda $ runs over the $n$-tuples of
partitions, that is non-increasing tuples on non-negative integers
$\lambda =(\lambda_1 \ge \lambda_2 \ge ...\ge \lambda _n ))$. The
coefficients $c_{\mu,\nu}^{\lambda}$ of the product
$$
s_{\mu } s_{\nu} =\sum c_{\mu ,\nu }^{\lambda}s_{\lambda },
$$
are known as the Littlewood-Richardson coefficients. These
coefficients occur in many problems, see, for example,
\cite{Fulton, M}.

From this definition obviously follow commutativity and
associativity
$$c_{\mu,\nu}^{\lambda}=c_{\nu,\mu}^{\lambda}
$$
\begin{equation}\label{ass1}
\sum_{\sigma}c_{\lambda,\mu}^{\sigma}c_{\sigma,\nu}^{\pi}=
\sum_{\tau}c_{\mu,\nu}^{\tau}c_{\lambda,\tau}^{\pi}.
\end{equation}
Littlewood and Richardson proposed combinatorial interpretation of
the coefficient $c_{\mu,\nu}^{\lambda}$ as the cardinality of
so-called LR-tableaux,  skew tableaux of shape $\lambda/\mu$ and
weight $\nu$ with the reverse lattice words reading (see, for
example, \cite{M}).

In \cite{BZ} it was found a nice interpretation of  LR-tableaux as
special functions on a triangle grid, see also
\cite{Buch,izv,stekl,surv,KT}. Due to this line of approach,
LR-coefficients count integer-valued discretely concave functions
on the triangle grid with the prescribed boundary values
(\cite{izv,stekl}) or hives (\cite{Buch,KT}).  In these terms the
commutativity and associativity relations are non-trivial
statements and give a hint for existing natural bijections among
corresponding sets.

In \cite{KTW2} is proposed a proof of this relation in the
language of hives, that is establishing a bijection between two
sets of  pairs of functions on the triangle grid. The second part
of the proof in not transparent for us\footnote{Recently in
\cite{HK} a sketch of another proof is proposed.}  and here we
propose a simple proof of this bijection. In fact, we prove a bit
stronger result of existence a polarized polymatroidal discretely
concave function on the tetrahedron grid with prescribed boundary
values on two adjoint faces of the tetrahedron (which, of course,
has to be discretely concave functions on the corresponding
triangular grids), for precise definitions see Section 3. As a
consequence of this result we get the bijection (due to the
polarization property, that is the {\em octahedron rule}). It is
interesting to point out that we recognized the octahedron rule as
a stronger version of the octahedron exchange axiom (\cite{DCU-2})
or a part of the local exchange axiom (\cite{DW92,MSh99}). The
local exchange axiom allows one of equivalent ways to construct
theory of polymatroidal discretely concave functions. Note that in
dimensions $d\ge 3$ there is no needs for higher dimensional
octahedron recursions to define polymatroidal discretely concave
functions in $\mathbb Z^d$.

As a consequence of our main result, we also get a bijection
$$
DC(\mu,\nu;\lambda)\cong DC(\nu,\mu;\lambda),
$$
which correspond to the commutativity
$c_{\mu,\nu}^{\lambda}=c_{\nu,\mu}^{\lambda}$.

In the language of $A$-type crystals of Kashiwara, a similar
associativity bijection can be obtained using twice the canonical
commutativity isomorphism $R:A\otimes B\to B\otimes A$
(\cite{bicrys}) in the arrays category of crystals, namely
applying $R:O_{\mu}\otimes (O_{\nu}\otimes O_{\lambda})\to
(O_{\nu}\otimes O_{\lambda})\otimes O_{\mu}$ and then
$R:O_{\tau}\otimes O_{\mu}\to O_{\mu}\times O_{\tau}$ for each
irreducible component $O_{\tau}$ of  $O_{\nu}\otimes O_{\lambda}$
we get an associativity bijection. Conjecturally, these bijections
have to coincide (looks alike Conjecture 3 in \cite{Pak}).

{\bf Acknowledgements}. We thank Igor Pak for discussions and
comments and directing us to the preprint \cite{HK}.

\section{Discretely concave functions on two-dimensional grid}

We let  $\Delta_n  =\Delta _n (O,X,Y)$ to denote the
two-dimensional grid of size $n$, that is the set of points
$(i,j)\in\mathbb Z^2$ such that  $i\ge 0$, $j\ge 0$ и $i+j\le n$.
Here we depict the grid of size $4$.

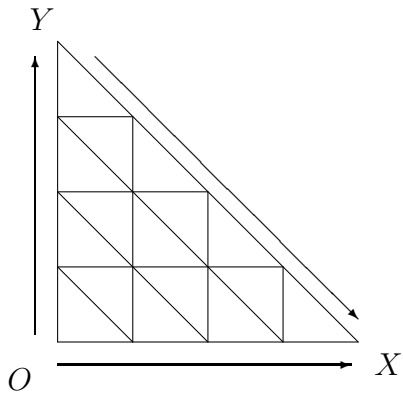
\unitlength=1mm \special{em:linewidth 0.4pt} \linethickness{0.4pt}
\begin{picture}(84.00,60.00)
\emline{40.00}{10.00}{1}{40.00}{50.00}{2}
\emline{40.00}{50.00}{3}{80.00}{10.00}{4}
\emline{80.00}{10.00}{5}{40.00}{10.00}{6}
\emline{50.00}{10.00}{7}{40.00}{20.00}{8}
\emline{40.00}{20.00}{9}{70.00}{20.00}{10}
\emline{70.00}{20.00}{11}{70.00}{10.00}{12}
\emline{70.00}{10.00}{13}{40.00}{40.00}{14}
\emline{40.00}{40.00}{15}{50.00}{40.00}{16}
\emline{50.00}{40.00}{17}{50.00}{10.00}{18}
\emline{60.00}{10.00}{19}{40.00}{30.00}{20}
\emline{40.00}{30.00}{21}{60.00}{30.00}{22}
\emline{60.00}{30.00}{23}{60.00}{10.00}{24}
\put(40.00,7.00){\vector(1,0){39.00}}
\put(37.00,11.00){\vector(0,1){37.00}}
\put(45.00,48.00){\vector(1,-1){35.00}}
\put(35.00,5.00){\makebox(0,0)[cc]{$O$}}
\put(38.00,53.00){\makebox(0,0)[cc]{$Y$}}
\put(84.00,7.00){\makebox(0,0)[cc]{$X$}}
\end{picture}

We consider a subclass of so called ``discretely concave''
functions defined on this grid $\Delta_n$ \cite{izv, stekl}.
Specifically, we cut the plane  $\mathbb R ^2$, which contains the
grid $\Delta_n$, by three types of lines $x=i$, $y=j$, $x+y=k$,
where $i,j$ and $k$ run over the integers. These lines cut the
triangle co$(\Delta _n )$ into small (unitary) triangles as we
depicted on the above Picture. Now, a given function $f$ defined
at the points of $\Delta_n$ we interpolate on each small triangle
by affinity. As a result of this interpolation we get a continuous
piece-wise linear function $\tilde f$ defined on the whole
triangle co$(\Delta _n )$.\medskip

{\bf Definition.} A function $f:\Delta _n\to\mathbb R$  is {\em
discretely concave} if the piece-wise linear interpolation $\tilde
f$ is a concave function\footnote{Of course, we can define
discretely concave functions on the whole lattice $f:\mathbb
Z^2\to\mathbb R\cup\{-\infty\}$, and due to the terminology in
\cite{DCU-2} such functions have to be called ``$\mathbb
A_2$-concave or two-dimensional polymatroidal concave function; in
\cite{MSh99} such functions called $M^{\natural}$-concave;
discretely concave functions on $\Delta_n$ called ``hives'' in
\cite{KT}. We prefer to use this terminology from \cite{stekl},
since on the one hand it reflects the discreteness of the domain
and on the another hand pointed out the crucial property of
``concavity''.}.\medskip

We can reformulate discrete concavity of a function $f$ without
using the interpolation $\tilde f$. Namely we have to require
validity of three types of ``rhombus'' inequalities

 (i)   $f(i,j)+f(i+1,j+1)\le f(i+1,j)+f(i,j+1)$;

(ii)  $f(i+1,j)+f(i+1,j+1)\ge f(i,j+1)+f(i+2,j)$;

(iii)  $f(i,j+1)+f(i+1,j+1)\ge f(i,j+2)+f(i+1,j)$.

That is for each small rhombus inside the grid (or the whole
$\mathbb Z^2$) the sum of values along the diagonal which belong
to a cut line is greater or equal to the sum along the diagonal
which does not belong.

Consider a function $f$ on the grid  $\Delta_n$ and consider its
restriction to each side of the triangle: the base of the
triangle, the left-hand side and the hypotenuse. Specifically, we
have to orient these sides as depict on the previous picture and
consider INCREMENTS of the function on each unit segment (arrow).
Then, increments along the  left-hand side constitute an $n$-tuple
$\mu =(\mu_1 ,\ldots,\mu _n )$. It is easy follows from the
rhombus inequalities of the type (i) and (iii) that  $\mu _1 \ge
\ldots\ge \mu _n$, that is $\mu$ is a partition. Analogously, we
consider the increments along the hypotenuse and will get a
non-increasing (the rhombus inequalities (ii) and (iii)) $n$-tuple
$\nu =(\nu _1 ,\ldots,\nu _n )$, and along the base of the
triangle we get non-increasing (the rhombus inequalities (i) and
(ii)) $n$-tuple и $\lambda =(\lambda_1 ,\ldots,\lambda _n )$. This
triple $\mu ,\nu $ and $\lambda $ we call boundary  {\em
increments} of a function $f$. Obviously, the increments are
invariant under adding a constant to $f$. Therefore, we have to
consider functions modulo adding a constant or to require
$f(0)=0$.

The set of all discretely concave functions on $\Delta_n$ with
increments  $(\mu ,\nu ,\lambda )$ we denote $DC_n (\mu ,\nu
,\lambda )$. This set is a polyhedron (possible empty) in the
space of all functions on  $\Delta _n$. Obvious necessary
conditions are non-increasing of  $\mu ,\nu $ and $\lambda $ and
validity of the equality  $\mu_1+...+\mu_n+ \nu
_1+...+\nu_n=\lambda _1+...+\lambda _n$. Of course for $n>2$ these
conditions are too far to be sufficient, for details see
\cite{Fulton,KT,stekl}.

For integers' tuples $(\mu ,\nu ,\lambda )$, of special interest
are integer points of this polytope, that is integer-valued
functions in $DC_n (\mu ,\nu ,\lambda )$. We let $DC_n (\mu ,\nu
,\lambda )(\mathbb Z)$ to denote this set. Within the array model
of $A$-type crystals \cite{bicrys}, elements of this set encode
the highest weight vectors in the tensor product of irreducible
crystals $O_{\mu}\otimes O_{\nu}$, which span crystals isomorphic
to $O_{\lambda}$. In particular, from that follows that
cardinality of this set coincides with the Littlewood-Richardson
coefficient $c_{\mu,\nu}^{\lambda}$, i.e.
\begin{equation}\label{lr}
|DC_n (\mu ,\nu ,\lambda )(\mathbb Z)|=c_{\mu,\nu}^{\lambda}.
\end{equation}
The first proof of  (\ref{lr}) is, implicitly, in \cite{BZ}, for
other proofs see Appendix in\cite{Buch}, \cite{izv,surv,KT}.

Thus, in the language of discretely concave functions,  the
associativity formula
$$
\sum_{\sigma}c_{\lambda,\mu}^{\sigma}c_{\sigma,\nu}^{\pi}=
\sum_{\tau}c_{\mu,\nu}^{\tau}c_{\lambda,\tau}^{\pi}
$$
reads as coincidence of cardinalities of the following two sets
$$
\coprod _{\lambda }DC_n (\mu ,\nu ,\lambda )(\mathbb Z )\times DC_n
(\lambda ,\pi ,\sigma )(\mathbb Z )
$$
and
$$
\coprod _{\tau }DC_n (\nu ,\pi ,\tau )(\mathbb Z )\times DC_n (\mu
,\tau ,\sigma )(\mathbb Z ).
$$
Of course, this gives a hint that it should exist a ``natural''
bijection between these sets. In \cite{KTW2} is proposed a
construction (we will call it KTW-construction) of a bijection
between sets
$$
\coprod _{\lambda }DC_n (\mu ,\nu ,\lambda )\times DC_n (\lambda
,\pi ,\sigma )
$$
and
$$\coprod _{\tau }DC_n(\nu ,\pi ,\tau )\times DC_n (\mu
,\tau ,\sigma )
$$
(here $\mu$, $\nu$, $\lambda $ etc. are real-valued, not necessary
integer-valued), which is, first, piece-wise linear and, second,
sends integer points into integer points\footnote{The final step
of the proof in \cite{KTW2} is not transparent and it is not clear
to us how Theorem 1.6 from \cite{FZ} can be used to finish the
proof.}. As a consequence of the main Theorem in the next section
we will get a proof of that the KTW-construction is a bijection
indeed.

\section{Functions on three dimensional grid}
Of course we have to recall nice and transparent KTW-construction.
For that we slightly transform the problem. Let we have two
functions  $f$ and  $f'$ on two copies of the grid $\Delta_n $
such that $\lambda $, the tuple of increments of $f$ along the
hypotenuse, is also the tuple of increments of $f'$ along the base
of the triangle. Then we can glue these  copies along the sides on
which $f$ and $f'$ have equal increments. As a result we get a
function on the ``square'' grid, which satisfies the rhombus
inequalities for almost all rhombus in this square grid, namely,
for all rhombus except those which have diagonals located on the
north-west south-east diagonal of square $\{0 \le x,y\le n\}$.
Now, given this function we have to construct a function on the
same square with the same boundary increments, but which is
piece-wise linear with respect to the mirror triangulation, see
the Picture.

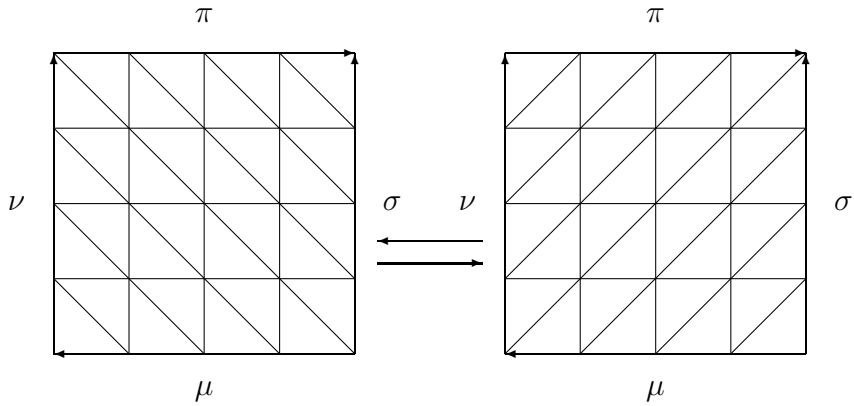
\unitlength=1mm \special{em:linewidth 0.4pt} \linethickness{0.4pt}
\begin{picture}(115.00,60.00)
\emline{20.00}{10.00}{1}{20.00}{50.00}{2}
\emline{10.00}{10.00}{3}{10.00}{50.00}{4}
\emline{10.00}{50.00}{5}{50.00}{50.00}{6}
\emline{50.00}{50.00}{7}{50.00}{10.00}{8}
\emline{50.00}{10.00}{9}{10.00}{10.00}{10}
\emline{30.00}{10.00}{11}{30.00}{50.00}{12}
\emline{40.00}{50.00}{13}{40.00}{10.00}{14}
\emline{50.00}{20.00}{15}{10.00}{20.00}{16}
\emline{10.00}{30.00}{17}{50.00}{30.00}{18}
\emline{50.00}{40.00}{19}{10.00}{40.00}{20}
\emline{80.00}{10.00}{21}{80.00}{50.00}{22}
\emline{70.00}{10.00}{23}{70.00}{50.00}{24}
\emline{70.00}{50.00}{25}{110.00}{50.00}{26}
\emline{110.00}{50.00}{27}{110.00}{10.00}{28}
\emline{110.00}{10.00}{29}{70.00}{10.00}{30}
\emline{90.00}{10.00}{31}{90.00}{50.00}{32}
\emline{100.00}{50.00}{33}{100.00}{10.00}{34}
\emline{110.00}{20.00}{35}{70.00}{20.00}{36}
\emline{70.00}{30.00}{37}{110.00}{30.00}{38}
\emline{110.00}{40.00}{39}{70.00}{40.00}{40}
\emline{20.00}{10.00}{41}{10.00}{20.00}{42}
\emline{10.00}{30.00}{43}{30.00}{10.00}{44}
\emline{40.00}{10.00}{45}{10.00}{40.00}{46}
\emline{10.00}{50.00}{47}{50.00}{10.00}{48}
\emline{20.00}{50.00}{49}{50.00}{20.00}{50}
\emline{50.00}{30.00}{51}{30.00}{50.00}{52}
\emline{40.00}{50.00}{53}{50.00}{40.00}{54}
\emline{70.00}{40.00}{55}{80.00}{50.00}{56}
\emline{90.00}{50.00}{57}{70.00}{30.00}{58}
\emline{70.00}{20.00}{59}{100.00}{50.00}{60}
\emline{110.00}{50.00}{61}{70.00}{10.00}{62}
\emline{80.00}{10.00}{63}{110.00}{40.00}{64}
\emline{110.00}{30.00}{65}{90.00}{10.00}{66}
\emline{100.00}{10.00}{67}{110.00}{20.00}{68}
\put(50.00,10.00){\vector(-1,0){40.00}}
\put(10.00,10.00){\vector(0,1){40.00}}
\put(10.00,50.00){\vector(1,0){40.00}}
\put(50.00,10.00){\vector(0,1){40.00}}
\put(110.00,10.00){\vector(-1,0){40.00}}
\put(70.00,10.00){\vector(0,1){40.00}}
\put(70.00,50.00){\vector(1,0){40.00}}
\put(110.00,10.00){\vector(0,1){40.00}}
\put(30.00,5.00){\makebox(0,0)[cc]{$\mu$}}
\put(5.00,30.00){\makebox(0,0)[cc]{$\nu$}}
\put(30.00,55.00){\makebox(0,0)[cc]{$\pi$}}
\put(55.00,30.00){\makebox(0,0)[cc]{$\sigma$}}
\put(90.00,5.00){\makebox(0,0)[cc]{$\mu$}}
\put(65.00,30.00){\makebox(0,0)[cc]{$\nu$}}
\put(90.00,55.00){\makebox(0,0)[cc]{$\pi$}}
\put(115.00,30.00){\makebox(0,0)[cc]{$\sigma$}}
\put(67.00,25.00){\vector(-1,0){14.00}}
\put(53.00,22.00){\vector(1,0){14.00}}
\end{picture}

Let us consider these $4$ triangles as four faces of the standard
tetrahedron  $\Delta _n (O,X,Y,Z)$ of size $n$.

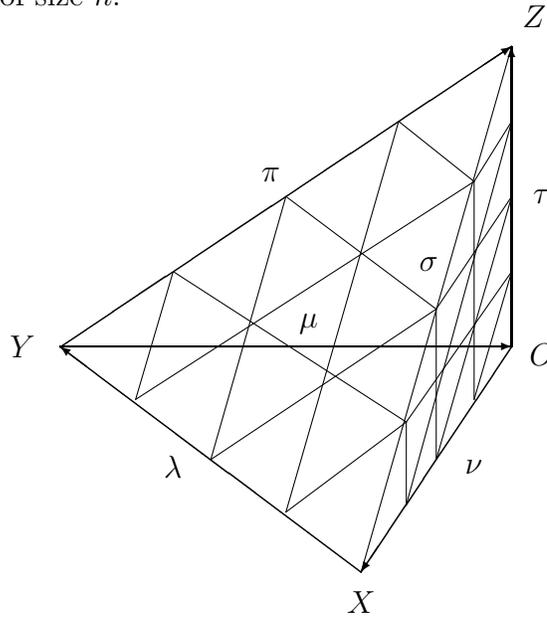
\unitlength=1mm \special{em:linewidth 0.4pt} \linethickness{0.4pt}
\begin{picture}(94.00,84.00)
\emline{30.00}{40.00}{1}{90.00}{40.00}{2}
\emline{90.00}{40.00}{3}{70.00}{10.00}{4}
\emline{70.00}{10.00}{5}{30.00}{40.00}{6}
\emline{30.00}{40.00}{7}{90.00}{80.00}{8}
\emline{90.00}{80.00}{9}{90.00}{40.00}{10}
\emline{70.00}{10.00}{11}{90.00}{80.00}{12}
\emline{90.00}{70.00}{13}{85.00}{62.00}{14}
\emline{85.00}{62.00}{15}{75.00}{70.00}{16}
\emline{90.00}{60.00}{17}{80.00}{45.00}{18}
\emline{80.00}{45.00}{19}{60.00}{60.00}{20}
\emline{90.00}{50.00}{21}{76.00}{30.00}{22}
\emline{76.00}{30.00}{23}{45.00}{50.00}{24}
\emline{40.00}{33.00}{25}{85.00}{62.00}{26}
\emline{85.00}{62.00}{27}{85.00}{33.00}{28}
\emline{50.00}{25.00}{29}{80.00}{45.00}{30}
\emline{80.00}{45.00}{31}{80.00}{25.00}{32}
\emline{80.00}{25.00}{33}{90.00}{60.00}{34}
\emline{60.00}{18.00}{35}{76.00}{30.00}{36}
\emline{76.00}{30.00}{37}{76.00}{19.00}{38}
\emline{76.00}{19.00}{39}{90.00}{70.00}{40}
\emline{90.00}{50.00}{41}{85.00}{33.00}{42}
\emline{75.00}{70.00}{43}{60.00}{18.00}{44}
\emline{50.00}{25.00}{45}{60.00}{60.00}{46}
\emline{45.00}{50.00}{47}{40.00}{33.00}{48}
\put(70.00,10.00){\vector(-4,3){40.00}}
\put(30.00,40.00){\vector(3,2){60.00}}
\put(90.00,40.00){\vector(-2,-3){20.00}}
\put(30.00,40.00){\vector(1,0){60.00}}
\put(90.00,40.00){\vector(0,1){40.00}}
\put(94.00,39.00){\makebox(0,0)[cc]{$O$}}
\put(70.00,6.00){\makebox(0,0)[cc]{$X$}}
\put(25.00,40.00){\makebox(0,0)[cc]{$Y$}}
\put(93.00,84.00){\makebox(0,0)[cc]{$Z$}}
\put(63.00,43.00){\makebox(0,0)[cc]{$\mu$}}
\put(85.00,24.00){\makebox(0,0)[cc]{$\nu$}}
\put(45.00,24.00){\makebox(0,0)[cc]{$\lambda$}}
\put(58.00,63.00){\makebox(0,0)[cc]{$\pi$}}
\put(79.00,51.00){\makebox(0,0)[cc]{$\sigma$}}
\put(94.00,60.00){\makebox(0,0)[cc]{$\tau$}}
\end{picture}

Initial functions  $f$ and $f'$ are given on the {\em ground}
$\Delta _n (O,X,Y)$ and on the {\em ceiling} $\Delta _n (X,Y,Z)$,
respectively. The outcome functions are defined on the {\em walls}
$\Delta _n (O,Y,Z)$ and $\Delta _n (O,X,Z)$. (On the picture we
triangulated the ground and one of the walls, analogously have to
be triangulated the ceiling and another wall.) The
KTW-construction propagates a function from the ground and ceiling
to all integer points of the tetrahedron
$$
           \Delta _n (O,X,Y,Z)=\{(x,y,z)\in \mathbb Z ^3,
x,y,z\ge 0, x+y+z\le n\}.
          $$
This propagation procedure runs through the so-called {\em
octahedron recursion}. It is clear from this recursion that it is
defined by piece-wise linear functions and is invertible. But that
is not clear that starting from discretely concave functions $f$
and $f'$ on the ground and ceiling we will end up with discretely
concave functions on the walls of our ``house''. As we pointed
out, the non-trivial step of the proof in \cite{KTW2} is made by
addressing to the ``tropical Laurent polynomials'' statement in
\cite{FZ}, and how this works is still unclear to us.

We are going  to identify the class of functions on the
three-dimensional grid, which are obtained due to this recursion,
and to prove that they form a subclass of the class of
three-dimensional discretely concave functions on the grid $\Delta
_n (O,X,Y,Z)$. Discrete concavity of such functions on the faces
of the tetrahedron holds true due to the definition. \medskip

Let us consider functions defined on the three-dimensional grid
$\Delta _n (O,X,Y,Z)$. (The convex hull of this grid we let to
denote by the same symbol and let to call it the {\em
tetrahedron}.)

We cut $\mathbb R^3$ (and the tetrahedron) by four series of
planes  $x=i$, $y=j$, $z=k$, $x+y+z=l$, where  $i$, $j$, $k$ and
$l$ run over integers (from  $0$ until $n$ for the tetrahedron).
As a result the tetrahedron is cut into small ``unitary'' pieces:
simplexes and octahedrons. (The matter has to be clear from the
following picture with  $n=2$.)

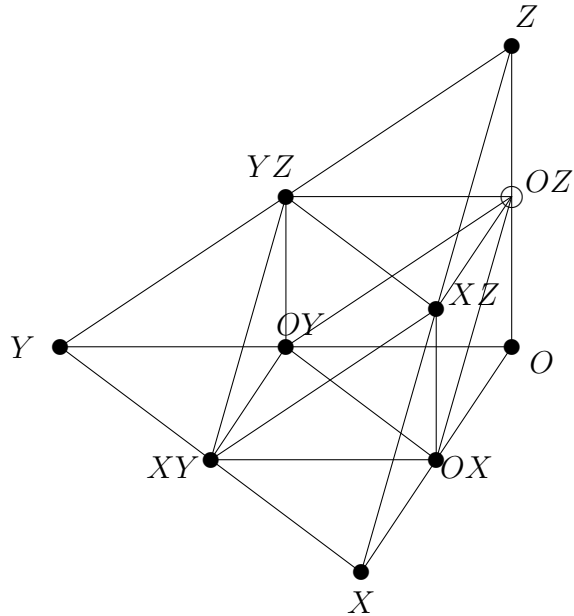
\unitlength=1mm \special{em:linewidth 0.4pt} \linethickness{0.4pt}
\begin{picture}(95.00,90.00)
\emline{30.00}{40.00}{1}{90.00}{40.00}{2}
\emline{90.00}{40.00}{3}{70.00}{10.00}{4}
\emline{70.00}{10.00}{5}{30.00}{40.00}{6}
\emline{30.00}{40.00}{7}{90.00}{80.00}{8}
\emline{90.00}{80.00}{9}{90.00}{40.00}{10}
\emline{70.00}{10.00}{11}{90.00}{80.00}{12}
\emline{90.00}{80.00}{13}{90.00}{80.00}{14}
\put(90.00,80.00){\circle*{2.00}} \put(60.00,60.00){\circle*{2.00}}
\put(30.00,40.00){\circle*{2.00}} \put(50.00,25.00){\circle*{2.00}}
\put(70.00,10.00){\circle*{2.00}} \put(80.00,45.00){\circle*{2.00}}
\put(90.00,40.00){\circle*{2.00}} \put(60.00,40.00){\circle*{2.00}}
\put(90.00,60.00){\circle{2.83}} \put(80.00,25.00){\circle*{2.00}}
\emline{80.00}{25.00}{15}{50.00}{25.00}{16}
\emline{50.00}{25.00}{17}{60.00}{40.00}{18}
\emline{60.00}{40.00}{19}{80.00}{25.00}{20}
\emline{80.00}{25.00}{21}{80.00}{45.00}{22}
\emline{80.00}{45.00}{23}{90.00}{60.00}{24}
\emline{90.00}{60.00}{25}{80.00}{25.00}{26}
\emline{50.00}{25.00}{27}{60.00}{60.00}{28}
\emline{60.00}{60.00}{29}{60.00}{40.00}{30}
\emline{60.00}{40.00}{31}{90.00}{60.00}{32}
\emline{90.00}{60.00}{33}{60.00}{60.00}{34}
\emline{60.00}{60.00}{35}{80.00}{45.00}{36}
\emline{80.00}{45.00}{37}{50.00}{25.00}{38}
\put(94.00,38.00){\makebox(0,0)[cc]{$O$}}
\put(70.00,6.00){\makebox(0,0)[cc]{$X$}}
\put(25.00,40.00){\makebox(0,0)[cc]{$Y$}}
\put(92.00,84.00){\makebox(0,0)[cc]{$Z$}}
\put(84.00,24.00){\makebox(0,0)[cc]{$OX$}}
\put(45.00,24.00){\makebox(0,0)[cc]{$XY$}}
\put(58.00,64.00){\makebox(0,0)[cc]{$YZ$}}
\put(95.00,62.00){\makebox(0,0)[cc]{$OZ$}}
\put(85.00,47.00){\makebox(0,0)[cc]{$XZ$}}
\put(62.00,43.00){\makebox(0,0)[cc]{$OY$}}
\end{picture}

The main attention we pay to the octahedrons. All octahedrons are
``similar'' and one is obtained from another by integer
translations. Each octahedron has three main diagonal, or
equivalently, three pairs of antipodal vertices: $OX$ and $YZ$,
$OY$ and $XZ$, $OZ$ and $XY$. The latter pair we distinguish among
others and the corresponding diagonal (parallel to the vector
$(1,1,-1)$) in each octahedron we call the {\em main diagonal}.

Reasons why we distinguish this  diagonal might be seen in the
simplest case  $n=2$. In this case, there is the sole point $OZ$
which does not belong to the ground face and ceiling face, and we
have to propagate a function to this point. The rule for the
propagation, called the {\em octahedron rule} in \cite{KTW2},
prescribes to set
$$
            f(OZ)=\max(f(OX)+f(YZ),f(OY)+f(XZ))-f(XY).
$$
In words this rule (the octahedron recursion) says: we have to set
the value at the ``free'' end of the main diagonal such that the
sum of values at the vertices of this main diagonal has to equal
the maximum over the sums with respect to other (non-main)
diagonals.

This leads us to the following\medskip

    {\bf  Definition}. A function $f:\mathbb Z ^3\to\mathbb
R\cup\{-\infty\}$) (or $f:\Delta_n (O,X,Y,Z)\to\mathbb R$ ) is
said to be  {\em polarized}, if, for any unitary octahedron, the
sum of values of $f$ at the vertexes of the main diagonal is equal
to the maximum of the sums of $f$ on vertexes on the other
diagonals.\medskip

It is obvious, that for any values at the points of the ground and
ceiling faces of the grid $\Delta_n (O,X,Y,Z)$, there exists a
unique polarized function on  $\Delta _n (O,X,Y,Z)$ with these
boundary values. This is the essence of the octahedron recursion.

We are going to prove that discrete concave ``initial'' data at
the ground face $OZY$ and the ceiling face  $XYZ$ produce a
polarized function which possesses the three-dimensional discrete
concavity. Without going deeply in details of this concept (which
has sense in any dimension, see \cite{DCU-2,MSh99}), we give
definitions which are of use here.

Recall that discrete concavity in dimension 2 is related to the
rhombus inequalities. In dimension 3, we have to require
inequalities for the octahedrons and the rhombus inequalities for
each of $4$ types cutting planes. Namely, in each such a cutting
plane we have three types of rhombus and the rhombus inequalities
require the sum at the vertices on the rhombus diagonal lying on a
cutting plane exceed the values at the another diagonal. The
octahedrons' inequalities says that values at two diagonal have to
be equal to the maximum over three values at the pairs of
antipodal vertexes. Functions with these properties are called
polymatroidal discretely concave functions in \cite{DCU-2} and
$M^{\natural}$-concave in \cite{MSh99}. The octahedron rule is a
stronger requirement and  implies validity of the octahedron
inequalities.
\medskip

{\bf Definition}. A function $f$ on the three-dimensional grid
$\Delta _n (O,X,Y,Z)$ is a {\em polarized discretely concave
function} if $f$ is polarized and all rhombus' inequalities hold
true. We let $PCPM_n$ to denote this set of functions.\medskip

It is clear that restrictions of a  $PCPM_n$-function  $f$ to any
cutting plane, and, hence, to the ground face, or to the ceiling
face, or to the walls, are two-dimensional discretely concave
functions. Thus, we will obtain that the KTW-construction is a
bijection if we get the following\medskip

     {\bf Theorem}. {\em Let $f$ be a polarized function on the grid
$\Delta_n  (O,X,Y,Z)$ and let the restrictions of $f$ to the
ground face and the ceiling face be two-dimensional discretely
concave functions. Then  $f\in PCPM_n$.}

\section{Proof of Theorem}

In the beginning we consider the case $n=2$ (see the picture
above). Let $f$ be a polarized function on $\Delta _2(O,X,Y,Z)$.
We have to check the rhombus inequalities. In this case all
rhombus are located on faces of the tetrahedron. For rhombus
located on the ground or ceiling faces, the corresponding
inequities hold true due to the assumptions. Thus, we have to
consider rhombus on the walls. Since the walls are symmetrical, it
suffices to verify the rhombus inequalities for the rhombus on the
wall $OXZ$. This wall contains three rhombus, the rhombus
$O,OZ,XZ,OX$, the rhombus $X,OX,OZ,XZ$ and the rhombus
$Z,OZ,OX,XZ$. Since the latter two rhombus are similar, we have to
verify the inequalities for the first and the second.

The first rhombus. We have to check that
$$
f(OZ)+f(OX)\ge f(XZ)+f(O).
$$
Since the vertexes  $OZ$ and  $XY$ form the main diagonal in the
octahedron, from the polarization condition, we have the
inequality
$$
f(OZ)+f(XY)\ge f(XZ)+f(OY).
$$
From the rhombus inequality, for the rhombus  $O,OX,XY,OY$ on the
ground face, we have
$$
f(OX)+f(OY)\ge f(O)+f(XY).
$$
Summing up these two, we get the required inequality.

The second rhombus. We have to show
$$
f(OX)+f(XZ)\ge f(X)+f(OZ).
$$
From the polarization we have either the equality
$$
f(OZ)+f(XY)=f(XZ)+f(OY),
$$
or
$$
f(OZ)+f(XY)=f(OX)+f(YZ).
$$
Assume the first equality holds true. Then, from the inequality
for the rhombus $X,XY,OY,OX$ lying on the ground face, we get the
inequality
$$
f(OY)+f(X)\le f(XY)+f(OX).
$$
Summing up it with the first equality, we get the required
inequality.

Assume the second equality holds true. Then using the rhombus
$X,XY,YZ,XZ$ from the ceiling face, we have the rhombus inequality
$$
f(X)+f(YZ)\le f(XY)+f(XZ).
$$
Summing up this inequality and the second equality, we get the
required inequality.

Thus, the case  $n=2$ is proven.\medskip

{\em The general case}. For this case, we need the following
\medskip

{\bf Lemma}. {\em  Let a polarized function  $f$ be such that its
restrictions to the ground and ceiling faces of $\Delta
_n(O,X,Y,Z)$ are discretely concave functions. Then its
restriction to any cutting plane parallel to the ground or ceiling
face, i.e. sets of the form either $\{z=a\}\cap \Delta
_n(O,X,Y,Z)$ or $\{x+y+z=d\}\cap \Delta _n(O,X,Y,Z)$, $a$,
$c\in\mathbb Z$, is a two-dimensional discrete concave fucntion.
}\medskip

{\em Proof}. We begin by checking the rhombus inequalities for
rhombus located in the plane $z=1$, and which have non-empty
intersection with the ceiling face. This case might be handled by
the case $n=3$, see the next picture.

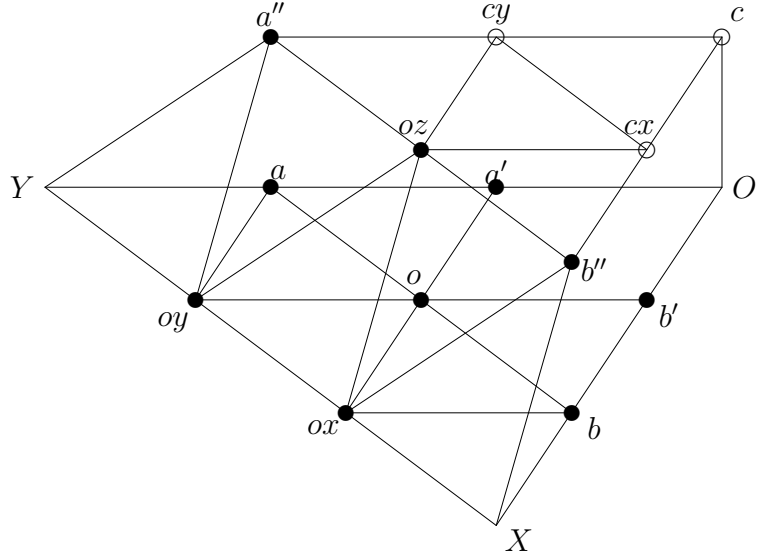
\unitlength=1mm \special{em:linewidth 0.4pt} \linethickness{0.4pt}
\begin{picture}(113.00,85.00)
\emline{20.00}{50.00}{1}{110.00}{50.00}{2}
\emline{110.00}{50.00}{3}{80.00}{5.00}{4}
\emline{80.00}{5.00}{5}{20.00}{50.00}{6}
\emline{20.00}{50.00}{7}{50.00}{70.00}{8}
\emline{50.00}{70.00}{9}{110.00}{70.00}{10}
\emline{110.00}{70.00}{11}{90.00}{40.00}{12}
\emline{90.00}{40.00}{13}{50.00}{70.00}{14}
\emline{110.00}{70.00}{15}{110.00}{50.00}{16}
\emline{90.00}{40.00}{17}{80.00}{5.00}{18}
\emline{80.00}{70.00}{19}{70.00}{55.00}{20}
\emline{70.00}{55.00}{21}{100.00}{55.00}{22}
\emline{100.00}{55.00}{23}{80.00}{70.00}{24}
\emline{50.00}{50.00}{25}{40.00}{35.00}{26}
\emline{40.00}{35.00}{27}{100.00}{35.00}{28}
\emline{90.00}{20.00}{29}{50.00}{50.00}{30}
\emline{90.00}{20.00}{31}{60.00}{20.00}{32}
\emline{60.00}{20.00}{33}{80.00}{50.00}{34}
\emline{50.00}{70.00}{35}{40.00}{35.00}{36}
\emline{40.00}{35.00}{37}{70.00}{55.00}{38}
\emline{70.00}{55.00}{39}{60.00}{20.00}{40}
\emline{60.00}{20.00}{41}{90.00}{40.00}{42}
\put(100.00,35.00){\circle*{2.00}} \put(90.00,20.00){\circle*{2.00}}
\put(60.00,20.00){\circle*{2.00}} \put(40.00,35.00){\circle*{2.00}}
\put(50.00,70.00){\circle*{2.00}} \put(70.00,55.00){\circle*{2.00}}
\put(90.00,40.00){\circle*{2.00}} \put(80.00,50.00){\circle*{2.00}}
\put(50.00,50.00){\circle*{2.00}} \put(80.00,70.00){\circle{2.00}}
\put(100.00,55.00){\circle{2.00}} \put(110.00,70.00){\circle{2.00}}
\put(70.00,35.00){\circle*{2.00}}
\put(113.00,50.00){\makebox(0,0)[cc]{$O$}}
\put(83.00,3.00){\makebox(0,0)[cc]{$X$}}
\put(17.00,50.00){\makebox(0,0)[cc]{$Y$}}
\put(112.00,73.00){\makebox(0,0)[cc]{$c$}}
\put(80.00,73.00){\makebox(0,0)[cc]{$cy$}}
\put(99.00,58.00){\makebox(0,0)[cc]{$cx$}}
\put(69.00,38.00){\makebox(0,0)[cc]{$o$}}
\put(37.00,32.00){\makebox(0,0)[cc]{$oy$}}
\put(57.00,18.00){\makebox(0,0)[cc]{$ox$}}
\put(69.00,58.00){\makebox(0,0)[cc]{$oz$}}
\put(51.00,52.00){\makebox(0,0)[cc]{$a$}}
\put(80.00,52.00){\makebox(0,0)[cc]{$a'$}}
\put(50.00,73.00){\makebox(0,0)[cc]{$a''$}}
\put(93.00,18.00){\makebox(0,0)[cc]{$b$}}
\put(103.00,33.00){\makebox(0,0)[cc]{$b'$}}
\put(93.00,39.00){\makebox(0,0)[cc]{$b''$}}
\end{picture}

On the first floor ($z=1$) there is two essentially different
rhombus: $c,cx,oz,cy$ and $a'',oz,cz,cy$. We have to check two
corresponding inequalities:
$$
f(cy)+f(cx)\ge f(c)+f(oz)
$$
and
$$
f(cy)+f(oz)\ge f(a'')+f(cx).
$$
The first inequality: Since $c$ is one of the vertexes of the main
diagonal $\{c,o\}$, the sum  $f(c)+f(o)$ is equal to either
$f(cy)+f(b')$ or $f(cx)+f(a')$. Suppose  $f(c)+f(o)=f(cx)+f(a')$
holds true. Since  $cy$ is the vertex of the main diagonal $\{cy,
oy\}$, we have  $f(cy)+f(oy)\ge f( a)+f(oz)$. Finally, for the
ground rhombus  $o,oy,a,a'$, we have the rhombus inequality
$f(a)+f(o)\ge f(a')+f(oy)$. Summing up these inequalities and the
equality $f(cx)+f(a')=f( c)+f(o)$, we get the desired inequality.
The case $f(c)+f(o)=f(cy)+f(b')$ is handled analogously.

The second inequality: Since the points  $cy$ and  $oy$ constitute
the main diagonal, we have, from the polarization, $f(cy)+f(oy)\ge
f(a'')+f(o)$. Now, for the main diagonal $\{cx, ox\}$, we have at
least one of two equalities either $f(cx)+f(ox)=$ is equal to
$f(b)+f(oz)$, or $f(cx)+f(ox)=f(o)+f(b'')$.

Suppose the first equality $f(b)+f(oz)= f(cx)+f(ox)$ holds true.
From the ground rhombus $\{ox,oy,o,b\}$ we have $f(o)+f(ox)\ge
f(oy)+f(b)$. Summing up  $f(cy)+f(oy)\ge f(a'')+f(o)$, $f(
b)+f(oz)= f(cx)+f(ox)$ and $f(o)+f(ox)\ge f(oy)+f(b)$, we get the
required inequality.

Suppose the second equality $f(o)+f(b'')=f(cx)+f(ox)$ holds true.
In this case, we use a rhombus inequality for rhombus $\{ox,ox,oz,
b''\}$ located on the CEILING face. Namely, summing up the
inequality  $f(ox)+f(oz)\ge f(oy)+f(b'')$ and the above equality
we get the required second inequality.

Now, one can see that rhombus inequality for the rhombuses  being
translations of the rhombus $\{oz,cy,c,cx\}$ by integer vectors,
follow from the polarization and rhombus inequalities on the
ground face. Thus, we get these inequalities. For a rhombus being
a translation of the rhombus $\{b'',oz,cy,cx\}$ or
$\{a'',cy,cx,oz\}$ we need rhombus inequalities for rhombuses
which have an edge on the ground face. Thus, this step allows us
to get rhombus inequalities for such rhombus which have a common
edge either on the ceiling face or the ground face. Thus step by
step, we obtain the rhombus inequalities for all rhombuses on the
planes parallel either the ground face or the ceiling face.
\medskip

Now all is prepared to prove the theorem. It remains to verify
rhombus inequalities for rhombuses located on planes parallel to
the walls. Pick a rhombus, say located on a plane parallel to the
wall $OXZ$. Then take the integer translation of the tetrahedron
$\Delta _2$ of size $2$, which contains this rhombus. By Lemma,
the function $f$ is discretely concave on the ceiling and ground
faces of this tetrahedron. Since $f$ is polarizes, from the proof
for the case $n=2$ follows that $f$ is a polarized discrete
concave on this tetrahedron and hence the desired rhombus
inequality holds true. \hfill Q.E.D.

\section{Commutativity}
Here we apply the main theorem to establish a bijective proof of
the commutativity
$$
c_{\mu,\nu}^{\lambda}=c_{\nu,\mu}^{\lambda}.
$$
Namely we obtain an isomorphism
\begin{equation}\label{comm1}
DC(\mu,\nu;\lambda)\cong DC(\nu,\mu;\lambda).
\end{equation}
Such kinds of isomorphisms in \cite{Pak} are called {\em
fundamental symmetries}. In \cite{bicrys}, this kind of symmetry
was obtained using the natural isomorphism $R:A\otimes B\to
B\otimes A$ in the array category of crystals. In \cite{HK} such
kind of a bijection is called commutor.

So, let $f\in DC(\mu,\nu;\lambda)$ be a discretely concave
function on the triangle grid $\Delta_n$. Consider the top half of
the octahedron inscribed in the tetrahedron
$\Delta_{2n}(O,X,Y,Z)$. This half-octahedron is a polymatroid and
there exist polymatroidal discretely concave functions on this
set\footnote{Recall that polymatroidal discretely concave
functions have efficiency domains of the form of polymatroids
(integer points), and a polymatroidal DC function being restricted
to a polymatroid, which is a subset of the efficiency domain,
remains polymatroidal DC function. We understand restriction as
setting $-\infty$ outside the restriction set.}. Let us call a
polymatroidal discretely concave function $f$ defined on a
polymatroidal set\footnote{A subset $P\subset \mathbb Z^3$ is a
polymatroidal set if $co (P)$ is a polymatroid, i.e. each edge of
$co(P)$ is parallel to some vector of the set $\{e_i, e_i-e_j\}$,
$i$, $j=1,2,3$, and $P=co (P)\cap\mathbb Z^3$.} $P\subset \mathbb
Z^3$ a {\em polarized} if the octahedron rule holds true for each
octahedron which has all vertices in $P$.

Now, we are looking for a polarized discretely concave function on
this half of the octahedron, which is equal to $f$ on the face of
this half of the octahedron which is located on the ceiling of
$\Delta_{2n}(O,X,Y,Z)$ and equals the separable function
$p_{\mu}\in DC(0,\mu;\mu)$ on the face $\{XY,YZ,OY\}$ (see the
next Picture).

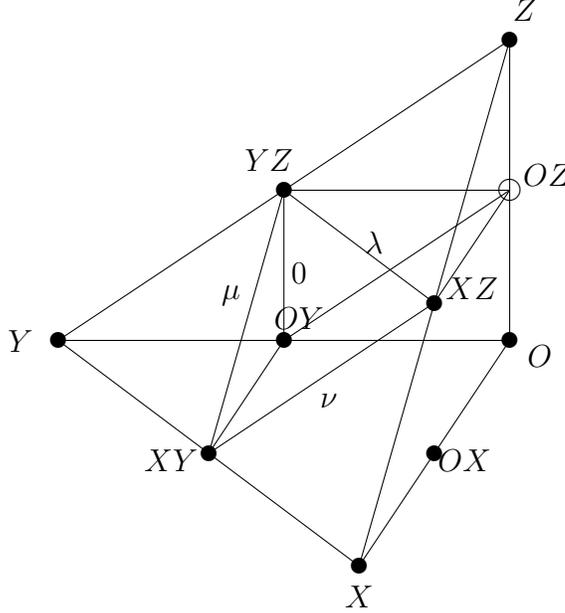
\unitlength=1.00mm \special{em:linewidth 0.4pt}
\linethickness{0.4pt}
\begin{picture}(95.00,84.00)
\emline{30.00}{40.00}{1}{90.00}{40.00}{2}
\emline{90.00}{40.00}{3}{70.00}{10.00}{4}
\emline{70.00}{10.00}{5}{30.00}{40.00}{6}
\emline{30.00}{40.00}{7}{90.00}{80.00}{8}
\emline{90.00}{80.00}{9}{90.00}{40.00}{10}
\emline{70.00}{10.00}{11}{90.00}{80.00}{12}
\emline{90.00}{80.00}{13}{90.00}{80.00}{14}
\put(90.00,80.00){\circle*{2.00}}
\put(60.00,60.00){\circle*{2.00}}
\put(30.00,40.00){\circle*{2.00}}
\put(50.00,25.00){\circle*{2.00}}
\put(70.00,10.00){\circle*{2.00}}
\put(80.00,45.00){\circle*{2.00}}
\put(90.00,40.00){\circle*{2.00}}
\put(60.00,40.00){\circle*{2.00}} \put(90.00,60.00){\circle{2.83}}
\put(80.00,25.00){\circle*{2.00}}
\emline{50.00}{25.00}{15}{60.00}{40.00}{16}
\emline{80.00}{45.00}{17}{90.00}{60.00}{18}
\emline{50.00}{25.00}{19}{60.00}{60.00}{20}
\emline{60.00}{60.00}{21}{60.00}{40.00}{22}
\emline{60.00}{40.00}{23}{90.00}{60.00}{24}
\emline{90.00}{60.00}{25}{60.00}{60.00}{26}
\emline{60.00}{60.00}{27}{80.00}{45.00}{28}
\emline{80.00}{45.00}{29}{50.00}{25.00}{30}
\put(94.00,38.00){\makebox(0,0)[cc]{$O$}}
\put(70.00,6.00){\makebox(0,0)[cc]{$X$}}
\put(25.00,40.00){\makebox(0,0)[cc]{$Y$}}
\put(92.00,84.00){\makebox(0,0)[cc]{$Z$}}
\put(84.00,24.00){\makebox(0,0)[cc]{$OX$}}
\put(45.00,24.00){\makebox(0,0)[cc]{$XY$}}
\put(58.00,64.00){\makebox(0,0)[cc]{$YZ$}}
\put(95.00,62.00){\makebox(0,0)[cc]{$OZ$}}
\put(85.00,47.00){\makebox(0,0)[cc]{$XZ$}}
\put(62.00,43.00){\makebox(0,0)[cc]{$OY$}}
\put(53.00,46.00){\makebox(0,0)[cc]{$\mu$}}
\put(66.00,32.00){\makebox(0,0)[cc]{$\nu$}}
\put(72.00,53.00){\makebox(0,0)[cc]{$\lambda$}}
\put(62.00,49.00){\makebox(0,0)[cc]{$0$}}
\end{picture}

Assume such a function $F$ exists, then since $F$ is a
polymatroidal discretely concave function, for the unitary
octahedrons which are cut by the face $\{XY,XY,OZ,OY\}$, the value
of $F$ at the vertex which does not belong  to the octahedron
$co\{XY,XY,OZ,OY, YZ\}$ equals $-\infty$. Thus for a rhombus on
the face $\{XY,XY,OZ,OY\}$ which obtained after cutting this face
by planes $z=a$, $x+y+z=b$, $a$ and $b\in\mathbb Z_+$, we have the
rhombus equality $F(A)+F(C)=F(B)+F(D)$ for the pair of antipodal
vertices $\{A,C\}$ and $\{B,D\}$ of the rhombus. From this follows
that increments of $F$ along  the edge $\{XZ,OZ\}$ is $\mu$ and
along the edge $\{OY,OZ\}$ is $\nu$. Since $F$ is the zero
constant on the edge $\{YZ,OY\}$, from discrete concavity on the
face $\{OY,YZ,OZ\}$ follows that the increments of $F$ along the
edge $\{YZ,OZ\}$ is $\nu$. Thus $F$ at the face $\{YZ,OZ,XZ\}$ is
indeed a function from $DC(\nu,\mu;\lambda)$. Since $F$ is
polarized we get a desired bijection.

{\bf Proposition}. {\em  There exists a polarized discretely
concave function on the half of the octahedron $co\{XY,XY,OZ,OY,
YZ\}$ such that its restriction to the face $\Delta_n(XZ,XY,XZ)$
is equal to a given function of $DC(\mu,\nu;\lambda)$ and its
restriction to $\Delta_n(OY,YZ,XY)$ is equal to the unique
function of $DC(0,\mu;\mu)$.}

{\em Proof}.  We proceed to find an appropriate polarized
discretely concave function on the whole simplex
$\Delta_{2n}(O,X,Y,Z)$. On the next two pictures we draw the
boundary values of such a function on the ceiling and ground
faces.

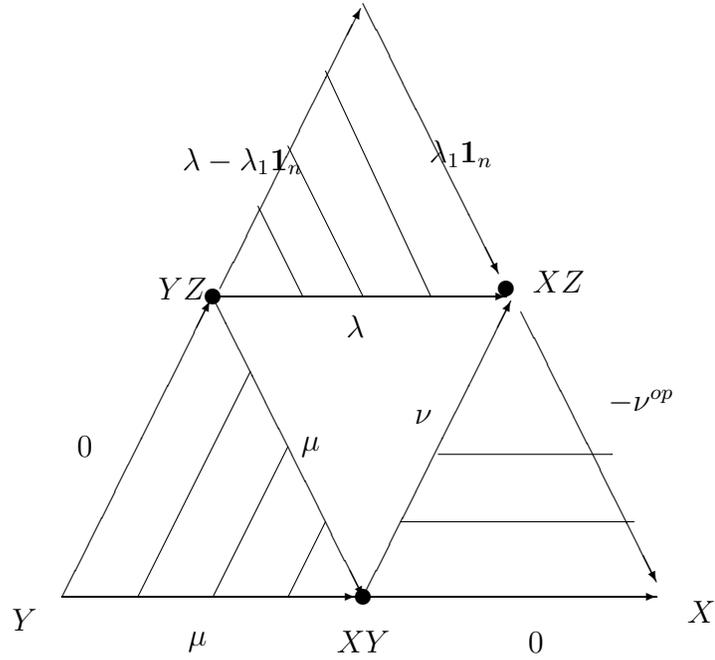
\unitlength=1mm \special{em:linewidth 0.4pt} \linethickness{0.4pt}
\begin{picture}(115.00,114.00)
\put(30.00,30.00){\vector(1,0){79.00}}
\put(50.00,70.00){\circle*{2.00}}
\put(89.00,71.00){\circle*{2.00}}
\put(70.00,30.00){\circle*{2.00}}
\put(50.00,70.00){\vector(1,-2){20.00}}
\put(70.00,30.00){\vector(1,2){19.67}}
\put(50.00,70.00){\vector(1,0){39.00}}
\put(30.00,30.00){\vector(1,0){39.00}}
\put(70.00,109.00){\vector(1,-2){18.00}}
\put(51.00,71.00){\vector(1,2){18.67}}
\put(91.00,68.00){\vector(1,-2){18.00}}
\put(25.00,27.00){\makebox(0,0)[cc]{$Y$}}
\put(67.00,114.00){\makebox(0,0)[cc]{$Z$}}
\put(115.00,28.00){\makebox(0,0)[cc]{$X$}}
\put(46.00,71.00){\makebox(0,0)[cc]{$YZ$}}
\put(96.00,72.00){\makebox(0,0)[cc]{$XZ$}}
\put(70.00,24.00){\makebox(0,0)[cc]{$XY$}}
\put(30.00,30.00){\vector(1,2){19.67}}
\put(48.00,24.00){\makebox(0,0)[cc]{$\mu$}}
\put(33.00,50.00){\makebox(0,0)[cc]{$0$}}
\put(63.00,50.00){\makebox(0,0)[cc]{$\mu$}}
\put(78.00,54.00){\makebox(0,0)[cc]{$\nu$}}
\put(69.00,66.00){\makebox(0,0)[cc]{$\lambda$}}
\put(54.00,88.00){\makebox(0,0)[cc]{$\lambda-\lambda_1{\mathbf
1}_n$}} \put(83.00,89.00){\makebox(0,0)[cc]{$\lambda_1{\mathbf
1}_n$}} \put(107.00,56.00){\makebox(0,0)[cc]{$-\nu^{op}$}}
\put(93.00,24.00){\makebox(0,0)[cc]{$0$}}
\emline{60.00}{30.00}{1}{65.00}{40.00}{2}
\emline{50.00}{30.00}{3}{60.00}{50.00}{4}
\emline{40.00}{30.00}{5}{55.00}{60.00}{6}
\emline{56.00}{82.00}{7}{62.00}{70.00}{8}
\emline{60.00}{90.00}{9}{70.00}{70.00}{10}
\emline{65.00}{100.00}{11}{79.00}{70.00}{12}
\put(84.00,58.00){\rule{12.00\unitlength}{0.00\unitlength}}
\emline{80.00}{49.00}{13}{103.00}{49.00}{14}
\emline{103.00}{49.00}{15}{103.00}{49.00}{16}
\emline{75.00}{40.00}{17}{106.00}{40.00}{18}
\end{picture}

On this Picture we draw an extension of a discretely concave
function defined on the grid $\Delta_n(YZ,XZ,XY)$ to the ceiling
$\Delta_{2n}(X,Y,Z)$. We draw lines outside the triangle $co
(YZ,XZ,XY)$ of constancy values of the extended function
($\nu^{op}$ denotes the non-decreasing tuple $(\nu_n,\ldots,
\nu_1)$ and ${\mathbf 1}_n:=(1,\ldots, 1)$).

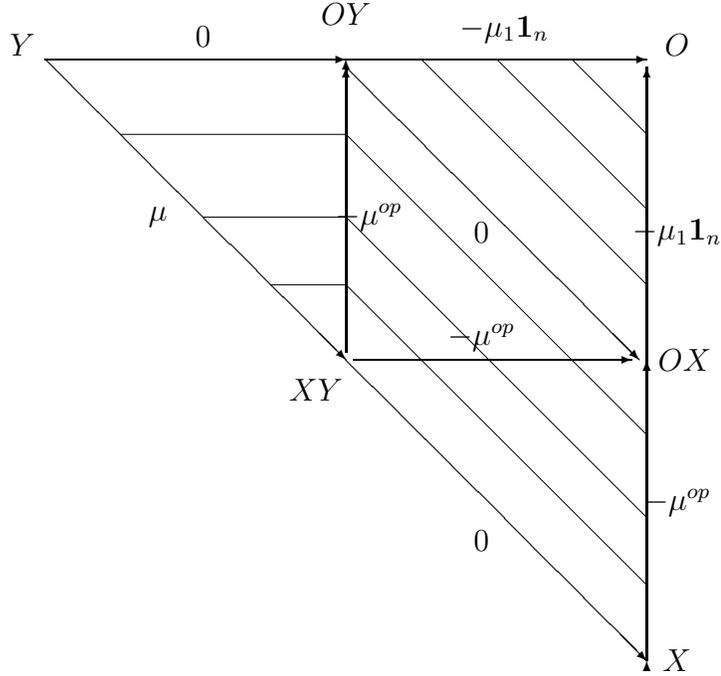
\unitlength=1mm \special{em:linewidth 0.4pt} \linethickness{0.4pt}
\begin{picture}(115.00,116.00)
\put(30.00,110.00){\vector(1,0){40.00}}
\put(70.00,110.00){\vector(0,0){0.00}}
\put(70.00,110.00){\vector(1,0){40.00}}
\put(30.00,110.00){\vector(1,-1){40.00}}
\put(70.00,70.00){\vector(1,-1){40.00}}
\put(110.00,30.00){\vector(0,0){0.00}}
\put(110.00,30.00){\vector(0,1){40.00}}
\put(110.00,70.00){\vector(0,0){0.00}}
\put(110.00,70.00){\vector(0,1){39.00}}
\put(70.00,71.00){\vector(0,1){38.00}}
\put(70.00,109.00){\vector(1,-1){39.00}}
\put(71.00,70.00){\vector(1,0){37.00}}
\put(27.00,112.00){\makebox(0,0)[cc]{$Y$}}
\put(114.00,112.00){\makebox(0,0)[cc]{$O$}}
\put(114.00,30.00){\makebox(0,0)[cc]{$X$}}
\put(70.00,116.00){\makebox(0,0)[cc]{$OY$}}
\put(115.00,70.00){\makebox(0,0)[cc]{$OX$}}
\put(66.00,66.00){\makebox(0,0)[cc]{$XY$}}
\put(45.00,89.00){\makebox(0,0)[cc]{$\mu$}}
\put(51.00,113.00){\makebox(0,0)[cc]{$0$}}
\put(73.00,89.00){\makebox(0,0)[cc]{$-\mu^{op}$}}
\put(88.00,73.00){\makebox(0,0)[cc]{$-\mu^{op}$}}
\put(88.00,46.00){\makebox(0,0)[cc]{$0$}}
\put(114.00,51.00){\makebox(0,0)[cc]{$-\mu^{op}$}}
\put(88.00,87.00){\makebox(0,0)[cc]{$0$}}
\put(114.00,87.00){\makebox(0,0)[cc]{$-\mu_1{\mathbf 1}_n$}}
\put(91.00,114.00){\makebox(0,0)[cc]{$-\mu_1{\mathbf 1}_n$}}
\emline{40.00}{100.00}{1}{70.00}{100.00}{2}
\emline{70.00}{100.00}{3}{110.00}{60.00}{4}
\emline{51.00}{89.00}{5}{70.00}{89.00}{6}
\emline{70.00}{89.00}{7}{110.00}{49.00}{8}
\emline{60.00}{80.00}{9}{70.00}{80.00}{10}
\emline{70.00}{80.00}{11}{110.00}{40.00}{12}
\emline{80.00}{110.00}{13}{110.00}{80.00}{14}
\emline{90.00}{110.00}{15}{110.00}{90.00}{16}
\emline{100.00}{110.00}{17}{110.00}{100.00}{18}
\end{picture}

On this Picture we draw the function on the ground face, the
boundary increments and the constancy levels determine a
discretely concave function on the ground face.

Thus, due to the main theorem, there exists a polarized discretely
concave function $f$ on $\Delta_{2n}(O,X,Y,Z)$, which has
restrictions to the ceiling and ground faces depicted on the above
two pictures.

Let us prove that this function $f$ restricted to the rhombus
$co\{X,XZ, OZ, OX\}$ (on the wall $\Delta_n(O,X,Z)$) is a
separable functions, i.e. for each unitary rhombus in this rhombus
and being homothetic to it, we have the rhombus equality. (Note
that these rhombuses are not located on the cutting planes.)

Assume this is done, then the tuple of increments along the side
$\{OY,OZ\}$ is $\nu$, hence that tuple along $\{YZ, OZ\}$ is also
$\nu$. Thus on the grid $\Delta_n(OZ,XZ,YZ)$ we get a function
from $DC (\nu,\mu;\lambda)$ and so the desired polarized function
on the half of the octahedron is obtained.

Thus, we have to prove validity of the rhombus equalities in the
rhombus $co\{X,XZ, OZ, OX\}$.

For this we use the following property of polarized functions on
$\Delta_2(O,X,Y,Z)$:

{\em Let $h\in PCPM_2$ and let $h(XY)+h(OX)=h(X)+h(OY)$. Then
$f(XZ)+f(OX)=f(X)+f(OZ)$ holds true}.

In fact, since there holds $h(XZ)+h(OX)\ge h(X)+h(OZ)$ and due to
$h(XY)+fh(OX)=h(X)+h(OY)$, we get $h(XY)+h(OZ)\le h(OY)+h(XZ)$.
Since $f\in PCPM_2$, the latter inequality holds as the reverse
inequality and so is the equality indeed, that is $f(XZ)+f(OX)=
f(X)+f(OZ)$ as claimed.

Now, having covered the trapezoid $co\{OY,O,X,XY\}$ by the ground
triangles of the all integer translations of $\Delta_2(O,X,Y,Z)$,
which intersect the interior of the trapezoid (see the next
Picture, where we draw some of screening triangles), and applying
the above claim to the corresponding tetrahedrons, we obtain the
same shape of the constancy levels of $f$ on $\{z=1\}\cap
co\{OY,O,X,XY, XZ,OZ\}$ as on the trapezoid $co\{OY,O,X,XY\}$.
Lifting higher and higher we get that $f$ satisfies all desired
rhombus equalities. Thus, this completes the proof of
Proposition.\bigskip

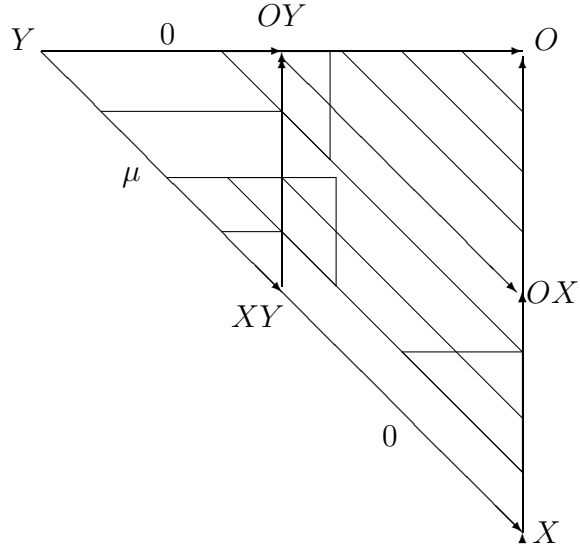
\unitlength=0.80mm \special{em:linewidth 0.4pt}
\linethickness{0.4pt}
\begin{picture}(115.00,116.00)
\put(30.00,110.00){\vector(1,0){40.00}}
\put(70.00,110.00){\vector(0,0){0.00}}
\put(70.00,110.00){\vector(1,0){40.00}}
\put(30.00,110.00){\vector(1,-1){40.00}}
\put(70.00,70.00){\vector(1,-1){40.00}}
\put(110.00,30.00){\vector(0,0){0.00}}
\put(110.00,30.00){\vector(0,1){40.00}}
\put(110.00,70.00){\vector(0,0){0.00}}
\put(110.00,70.00){\vector(0,1){39.00}}
\put(70.00,71.00){\vector(0,1){38.00}}
\put(70.00,109.00){\vector(1,-1){39.00}}
\put(27.00,112.00){\makebox(0,0)[cc]{$Y$}}
\put(114.00,112.00){\makebox(0,0)[cc]{$O$}}
\put(114.00,30.00){\makebox(0,0)[cc]{$X$}}
\put(70.00,116.00){\makebox(0,0)[cc]{$OY$}}
\put(115.00,70.00){\makebox(0,0)[cc]{$OX$}}
\put(66.00,66.00){\makebox(0,0)[cc]{$XY$}}
\put(45.00,89.00){\makebox(0,0)[cc]{$\mu$}}
\put(51.00,113.00){\makebox(0,0)[cc]{$0$}}
\put(88.00,46.00){\makebox(0,0)[cc]{$0$}}
\emline{40.00}{100.00}{1}{70.00}{100.00}{2}
\emline{70.00}{100.00}{3}{110.00}{60.00}{4}
\emline{51.00}{89.00}{5}{70.00}{89.00}{6}
\emline{70.00}{89.00}{7}{110.00}{49.00}{8}
\emline{60.00}{80.00}{9}{70.00}{80.00}{10}
\emline{70.00}{80.00}{11}{110.00}{40.00}{12}
\emline{80.00}{110.00}{13}{110.00}{80.00}{14}
\emline{90.00}{110.00}{15}{110.00}{90.00}{16}
\emline{100.00}{110.00}{17}{110.00}{100.00}{18}
\emline{79.00}{71.00}{19}{79.00}{89.00}{20}
\emline{79.00}{89.00}{21}{61.00}{89.00}{22}
\emline{61.00}{89.00}{23}{79.00}{71.00}{24}
\emline{110.00}{40.00}{25}{110.00}{60.00}{26}
\emline{110.00}{60.00}{27}{90.00}{60.00}{28}
\emline{90.00}{60.00}{29}{110.00}{40.00}{30}
\emline{60.00}{110.00}{31}{78.00}{92.00}{32}
\emline{78.00}{92.00}{33}{78.00}{110.00}{34}
\emline{78.00}{110.00}{35}{60.00}{110.00}{36}
\end{picture}

{\bf Remark}. It is not difficult to show that the bijection
constructed above coincides with the bijection in \cite{HK}.


\begin{thebibliography}{99}
\bibitem{BZ} Berenstein A. and A.Zelevinsky, Involutions on the
Gel'fand-Tsetlin patterns and muliplicities in the skew
$GL_n$-modules, {\em Doklady AN SSSR} (in Russian), 300(1988),
1291--1294, English translation in Soviet Math. Dokl. 37 (1988),
799--802
\bibitem{Buch} Buch A., The saturation conjecture (after A.Knutson
and T.Tao). With an appendix by William Fulton. {\em Enseign.
Math.} (2) 46 (2000), 43-60
\bibitem{izv}  Danilov V.I. and G. A. Koshevoy,
Discrete convexity and Nilpotent operators. {\it Math. Izvestiya
Russ. Acad. Sciences}, 67 (2003), 3--20 (in Russian) (English
transl. in Russian Acad. Sci. Izv. Math. 67:1 (2003))
\bibitem{stekl}  Danilov V.I. and G. A. Koshevoy,
Discrete convexity and Hermitian matrices. Proc. V.A.Steklov Inst.
Math.; Vol 241 (2003), 68--90 (English transl. in Number Theory,
Algebra and Algebraic Geometry: Collected papers. Dedicated to the
80th birthday of academician Igor' Rostislavovich Shafarevich.
Moscow. Nauka Publ. Proc. V.A.Steklov Inst. Math.; Vol 241,
2003)
\bibitem{DCU-2} Danilov V.I. and G.A. Koshevoy, Discrete convexity and
unimodularity II, (in preparation)
\bibitem{bicrys} Danilov V.I. and G.A. Koshevoy, Bicrystals and
crystal $GL(V), GL(W)$-duality. preprint RIMS-1485 (2004)
\bibitem{surv} Danilov V.I. and G.A. Koshevoy,
Arrays and combinatorics of Young tableaux. {\em
Uspehi Math. Nayk} (in Russian) (to appear)
\bibitem{DW92} Dress A.M.W. and W. Wenzel, Valuated Matroids,
{\em Advances in Mathematics}, 91 (1992), 158--208
\bibitem{FZ} Fomin S. and A.Zelevinsky, The Laurent phenomenon.
{\em Adv. in Appl. Math.}, 28 (2002), 119--144
\bibitem{Fulton} Fulton, W., Eigenvalues, Invariant Factors, Highest Weights,
and  Schubert Calculus, {\em Bull. Am. Math. Soc.,} 2000, vol. 37,
pp. 209--249.
\bibitem{HK} Henriques A. and J. Kamnitzer, The
octahedron recurrence and $gl_n$ ctystals. arXiv:math.CO/0408114 
\bibitem{KT}
Knutson, A. and Tao, T., The Honeycomb Model of $GL_n(\mathbb C)$ Tensor Product
 I: Proof of the Saturation Conjecture, {\em J. Am. Math. Soc.,} 1999, vol. 12,
 pp. 1055--1090.
\bibitem{KTW2} Knutson A., T.Tao and C. Woodward, A positive proof
of the Littlewood-Richardson rule using the octahedron reccurence.
arXiv:math.CO/0306274
\bibitem{M} { Macdonald I.G.} {\em Symmetric Functions and Hall
Polinomials}. Oxford Mathematical Monographs. Clarendon Press.
Oxford 1979
\bibitem {MSh99} Murota K., Discrete convex analysis,
{\em Mathematical Programming}, 83 (1998), 313--371
(see also Murota K., {\em Discrete convex analysis}, SIAM Monographs, 2003,
Philadelphia)
\bibitem{Pak} Pak.I and E.Vallejo, Reductions
of Young tableaux bijections. arXiv:math.CO/0408171

\end{thebibliography}
\end{document}